\documentclass[11pt]{amsproc} 

\usepackage{amsmath, amssymb, amsthm, setspace, mathtools, graphicx, enumitem, nicefrac}
\usepackage{thmtools}
\usepackage[usenames]{color}
\usepackage{tikz-cd}

\usepackage[pagebackref,colorlinks,linkcolor=blue,citecolor=blue,urlcolor=blue,hypertexnames=true]{hyperref}
\usepackage{amsrefs} 

\declaretheorem[style=plain,numberwithin=section,name=Theorem]{theorem}
\declaretheorem[style=plain,numberlike=theorem,name=Proposition]{proposition}
\declaretheorem[style=plain,numberlike=theorem,name=Lemma]{lemma}
\declaretheorem[style=plain,numberlike=theorem,name=Corollary]{corollary}

\declaretheorem[style=definition,numberlike=theorem,name=Question]{question}

\declaretheorem[style=remark,numberlike=theorem,name=Remark]{remark}
\declaretheorem[style=remark,numberlike=theorem,name=Example]{example}

\numberwithin{equation}{section}

\newcommand{\N}{\mathbb{N}}             
\newcommand{\Z}{\mathbb{Z}}             
\newcommand{\R}{\mathbb{R}}             
\newcommand{\SL}{\mathrm{SL}}           
\newcommand{\GL}{\mathrm{GL}}           
\newcommand{\St}{\mathrm{St}}           
\newcommand{\Aut}{\mathrm{Aut}}         
\newcommand{\SAut}{\mathrm{SAut}}       
\newcommand{\one}{\mathbf{1}}           
\newcommand{\transp}{\mathsf{T}}        
\newcommand{\ev}{\mathit{ev}}           
\newcommand{\id}{\mathit{id}}           
\newcommand{\I}{\mathrm{I}}             

\DeclareMathOperator{\Ima}{im}
\DeclareMathOperator{\tr}{tr}
\DeclareMathOperator{\Span}{Span}
\DeclareMathOperator{\interior}{int}

\author{Martin Nitsche}
\address{Martin Nitsche, Karlsruhe Institute of Technology, Germany}
\email{martin.nitsche@kit.edu}

\title{Computer proofs for Property~(T), and SDP duality}

\begin{document}

\onehalfspace

\begin{abstract}
We show that the semidefinite programs involved in the computer proofs for Kazhdan's property~(T) satisfy strong duality and that the dual programs have a geometric interpretation in terms of harmonic cocycles.
By dualizing geometric arguments about cocycles, we are able to simplify the property~(T) SDP in the case where it carries a symmetry by finite-order inner automorphisms.
As an application, we simplify the SDP proof for $\SL(n,\Z)$ and we prove that $\Aut(F_4)$ has property~(T).
\end{abstract}

\maketitle

\section{Introduction}
Kazhdan's property~(T) is a strong rigidity property for groups.
It has long been studied, and there are various equivalent characterizations of this property, highlighting different aspects (see \cite{BHV} for a textbook reference).
According to one common definition, a discrete group $\Gamma$, generated by a finite symmetric set $S=S^{-1}\subset\Gamma$, has property~(T), iff the Laplace operator in the maximal $\mathrm{C}^*$-algebra, $\Delta=\Delta_S\vcentcolon=|S|\cdot\one-\sum_{s\in S} s\in\mathrm{C}^*_\mathrm{max}(\Gamma)$, has a spectral gap directly above zero.

Traditionally, it is difficult to prove that a given group satisfies property~(T).
But in recent years Ozawa's article \cite{O} has kicked off a new approach to proving property~(T) with the computer.
Ozawa showed that if the Laplacian has a spectral gap, this fact is witnessed in the group algebra $\R[\Gamma]$.
Namely, there exist $\varepsilon>0\in\R$ and $\omega_1,\dots,\omega_n$ in the augmentation ideal $\I[\Gamma]\vcentcolon=\big\{\sum_i\lambda_i\gamma_i\,\big\vert\,\gamma_i\in\Gamma,\lambda_i\in\R,\sum_i\lambda_i=0\big\}$, such that
\begin{equation}
\Delta^2-\varepsilon\Delta=\sum_{i=1}^n \omega_i^*\omega_i,\label{eqn:ozawa-witness}
\end{equation}
where the involution $\cdot^*$ is given by $\sum_i\lambda_i\gamma_i\mapsto\sum_i\lambda_i\gamma_i^{-1}$.
Since elements in $\R[\Gamma]$ have finite support, the search for a property~(T) witness $\omega_1,\dots,\omega_n$ can be done by the computer.

The first implementation of a property~(T) computer proof was worked out for $\SL(3,\Z)$ by Netzer and Thom \cite{NT} in the framework of semidefinite programs (SDP).
Subsequently, the SDP approach has been used to both prove property~(T) \cites{KNO,KKN,CCKW} and to estimate the spectral gap size of various property~(T) group Laplacians \cites{F,KN}.
Although, in theory, the computer proof will eventually succeed for any property~(T) group, the challenge is to make the computation feasible in practice.

\medskip

The SDP proofs for property~(T) are remarkable in both their simplicity and their success in practice.
One downside, however, is that the computer proofs do not provide much human insight into why the group in question has property~(T), since typically there does not appear any humanly recognizable pattern among the coefficients of the witness that the computer produces.
In this article we take a first step in the direction of making the property~(T) computer proofs more accessible to humans by approaching the property~(T) SDP from another, more geometric perspective.

By Shalom's Theorem a finitely generated group $\Gamma$ has property~(T) exactly if the reduced cohomology $\bar{H}^1(\Gamma,\pi)$ is trivial for any orthogonal representation $\pi$ of $\Gamma$ on a Hilbert space $\mathcal{H}$. As we will recall in Section~\ref{sec:cocycle-recap}, this can be interpreted geometrically as saying that for no representation $\pi$ does there exist a non-trivial harmonic cocycle $c\colon\Gamma\to\mathcal{H}$. Combining Shalom's Theorem with Ozawa's result, one hence obtains the following dichotomy:

\begin{proposition}\label{prp:dichotomy}
Let $\Gamma$ be a finitely generated group and $S=S^{-1}$ a symmetric finite generating set. Then either
\begin{enumerate}
\item
there exists a witness for the spectral gap of the Laplacian, as in Equation~\ref{eqn:ozawa-witness}, or
\item
there exists an orthogonal representation and a non-trivial 1-cocycle that is harmonic, i.e. that satisfies $\sum_{s\in S}c(s)=0$.
\end{enumerate}
\end{proposition}
The property~(T) SDP, as introduced by Ozawa and Netzer--Thom, was motivated by the search for a spectral gap witness.
But as we will see, the search for a harmonic cocycle also leads to an SDP, which is exactly the dual of the spectral gap SDP.
Computer programs that solve the spectral gap SDP will typically also solve dual cocycle SDP, at the same time maximizing the size of the spectral gap that is witnessed and minimizing the curvature of a (partial) cocycle.

In Section~\ref{sec:T-duality} we show that the spectral gap SDP and the cocycle SDP satisfy what is called strong duality. This reproves Proposition~\ref{prp:dichotomy} in a direct and fairly elementary way. By combining strong duality with Ozawa's result \cite{O}*{Prop.\ 4} we also obtain a quantitative version of Proposition~\ref{prp:dichotomy}.

\newtheorem*{restate:quantitative}{Theorem~\ref{thm:quantitative-version}}
\begin{restate:quantitative}
Let $\Gamma$ be a finitely generated group, $S=S^{-1}$ a finite generating set and $\Delta\in\R[\Gamma]$ the corresponding Laplacian.
Let $\mathcal{C}$ denote the set of 1-cocycles $\Gamma\to\mathcal{H}$ satisfying $\sum_{s\in S}\|c(s)\|^2=2$ for any orthogonal representation $\Gamma\curvearrowright\mathcal{H}$.

Then the infimum $\inf_{c\in\mathcal{C}}\|\sum_{s\in S}c(s)\|^2$ is assumed and agrees with the size of the spectral gap of $\Delta$ above $0$.
\end{restate:quantitative}

The cocycle SDP leads to the exact same computer calculation as the spectral gap SDP. But the geometrical nature of the cocycle perspective opens up the calculation to geometric intuition. In Section~\ref{sec:dualizing-arguments} we demonstrate how geometric arguments can be used to simplify the SDP, making the calculation faster. Our simplification works in the case where $\Gamma$ contains a finite subgroup $H$ that normalizes the generating set $S$. We then examine the SDP property (T) proof for $\SL(n,\Z)$ from the geometric perspective and draw a connection to norm estimates for Harper's operator.

\medskip

The most important application of the SDP property~(T) proofs, to date, are the automorphism groups of the free groups, $\Aut(F_n)$. These groups had, for $n\geq 4$, long been suspected to satisfy property~(T) and the problem of proving this had been an important open question.
By exploiting symmetry, Kaluba--Nowak--Ozawa \cite{KNO} were able to simplify the SDP for $\Aut(F_5)$ far enough that the computer proof could be carried out in practice. By exploiting further structure of the SDP, Kaluba--Kielak--Nowak \cite{KKN} were able to extend the proof to $\Aut(F_n)$, $n\geq 6$. For $n=4$ they did not find a spectral gap witness within the search depth that was feasible to handle with their SDP implementation.

By using our simplification, we were able to efficiently cover a bigger search space in the search for a spectral gap witness, and by doing so we were able to prove property~(T) for the group $\Aut(F_4)$. This was the last $\Aut(F_n)$ group for which a proof had still been missing, as it is known that $\Aut(F_n)$ does not satisfy property~(T) for $n\in\{2,3\}$ (see the introduction to~\cite{KNO}).

\newtheorem*{restate:autF4}{Theorem~\ref{thm:autF4}}
\begin{restate:autF4}
$\Aut(F_4)$, the automorphism group of the free group over four generators, satisfies property~(T).
\end{restate:autF4}

\section{The cocycle perspective}\label{sec:cocycle-perspective}

\subsection{Preliminaries on cocycles and functions conditionally of positive type}\label{sec:cocycle-recap}
We begin by recalling the notion of 1-cocycles and the associated functions conditionally of positive type. We work with real Hilbert spaces and orthogonal representations, but one can easily pass between this setting and unitary representations on complex Hilbert spaces by complexifying or forgetting the complex structure. For a textbook account see~\cite{BHV}.

Let $\pi$ be an orthogonal representation of $\Gamma$ on a real Hilbert space $\mathcal{H}$. A 1-cocycle for $\pi$ is a map $c\colon\Gamma\to\mathcal{H}$ that sends the neutral element $\one$ to $0$ and satisfies the relation $c(\gamma_1\gamma_2)=c(\gamma_1)+\pi(\gamma_1)c(\gamma_2)$.
The set of 1-cocycles forms a vector space $Z^1$. Since every cocycle is determined by its values on $S$, we can think of $Z^1$ as a subspace of the Hilbert space $\mathcal{H}^{|S|}$.
The space of cocycles contains the closure of the image of the coboundary map $\delta^0\vcentcolon=\bigoplus_{s\in S}\id-\pi(s)\colon\mathcal{H}\to\mathcal{H}^{|S|}$. The quotient $Z^1/\overline{\Ima\delta^0}$ is denoted $\bar{H}^1(\Gamma,\pi)$ and can be identified with the space of \emph{harmonic cocycles}, $\bar{H}^1(\Gamma,\pi)\cong Z^1\cap(\Ima\delta^0)^\perp=Z^1\cap\ker(\delta^0)^*$.

Every 1-cocycle $c$ gives rise to an associated function $\bar{x}'\colon\Gamma\to\R$, $\gamma\mapsto -1/2\cdot\|c(\gamma)\|^2$ that maps the neutral element to zero.
For us it will be more convenient to work with the associated (unbounded) functional $x'\colon\I[\Gamma]\to\R$, $\gamma-\one\mapsto\bar{x}'(\gamma)$, which carries exactly the same information as $\bar{x}'$.
Let now $\Phi\colon\I[\Gamma]\otimes\I[\Gamma]\to\I[\Gamma]$ be the linear map that sends $\gamma_1\otimes\gamma_2$ to $\gamma_1^{-1}\gamma_2$.
The associated function $\bar{x}'$ of a 1-cocycle is always \emph{conditionally of positive type}, meaning that the bilinear form defined on $\I[\Gamma]$ by $\langle\cdot,\cdot\rangle_{x'}\vcentcolon=x'\circ\Phi$ is positive semidefinite. The form $\langle\cdot,\cdot\rangle_{x'}$ is automatically invariant under the left multiplication action of $\Gamma$ on $\I[\Gamma]$.

In the other direction, any $\Gamma$-invariant positive semidefinite bilinear form $\langle\cdot,\cdot\rangle_{x'}$ on $\I[\Gamma]$ factors as the concatenation of $\Phi$ with a function conditionally of positive type. For $\langle\cdot,\cdot\rangle_{x'}$ given, let $\mathcal{H}_{x'}$ denote the Hilbert space completion of $\I[\Gamma]/\{\|\cdot\|_{x'}=0\}$. The left multiplication action of $\Gamma$ on $\I[\Gamma]$ induces a $\Gamma$-action on $\mathcal{H}_{x'}$. The map $c_{x'}\colon\Gamma\to\mathcal{H}_{x'}$ that sends $\gamma\in\Gamma$ to the equivalence class of $\one-\gamma$ is a cocycle with respect to this action, and it satisfies $-1/2\cdot\|c_{x'}(\gamma)\|^2=\bar{x}'(\gamma)$.

The above procedure provides essentially a one-to-one correspondence between cocycles and $\Gamma$-invariant positive semidefinite forms on $\I[\Gamma]$.
Indeed, any cocycle with associated functional $x'$ is related to $c_{x'}\colon\Gamma\to\mathcal{H}_{x'}$ by an orthogonal equivalence of Hilbert spaces, after possibly enlarging the Hilbert space $\mathcal{H}_{x'}$. For details see \cite{BHV}*{proof of Thm.\ 2.10.2}.

\subsection{Proving property~(T) with partial cocycles}\label{sec:cocycle-sdp}

Shalom's theorem states that a finitely generated group $\Gamma$ has property~(T) exactly if $\bar{H}^1(\Gamma,\pi)=0$ for all orthogonal representations $\pi$.
Hence, to prove that $\Gamma$ has property~(T), one has to show that for no orthogonal representation $\pi$ does there exist a cocycle $c\colon\Gamma\to\mathcal{H}$ that is both harmonic and non-trivial.

Viewing $c$ as an element in $\mathcal{H}^{|S|}$, the harmonicity condition can be rewritten as
\begin{align*}
0=(\delta^0)^*(c)&=\sum_{s\in S}c(s)-\pi(s)^*c(s)=\sum_{s\in S}c(s)+c(s^{-1})-\big(c(s^{-1})+\pi(s^{-1})c(s)\big)\\&=\sum_{s\in S}c(s)+c(s^{-1})-c(s^{-1}s)=\sum_{s\in S}2c(s)-c(\one)=2\sum_{s\in S}c(s),
\end{align*}
using the fact that $S$ is symmetric.
As recapitulated above, we may assume that $c$ comes from a function $\bar{x}'$ conditionally of positive type, i.e.\ $c=c_{x'}\colon\Gamma\to\mathcal{H}_{x'}$.
In this setting the harmonicity condition becomes
\begin{equation}\label{eqn:harmonicity-condition}
0
=\|{\textstyle\sum}c(s)\|^2
=\|{\textstyle\sum}c_{x'}(\one)-c_{x'}(s)\|^2
=\|{\textstyle\sum}\one-s\|_{x'}^2
=\|\Delta\|_{x'}^2=x'(\Delta^2).
\end{equation}
We will call $\kappa(c)\vcentcolon=x'(\Delta^2)$ the \emph{curvature} of $c=c_{x'}$.

To enforce non-triviality, we will require that $c$ must be \emph{normalized} in the sense that
\begin{equation}\label{eqn:normalizing-condition}
2=\|c\|^2_{\mathcal{H}^{|S|}}=\sum_{s\in S}\|c(s)\|^2=\sum_{s\in S}\|c_{x'}(\one)-c_{x'}(s)\|^2=\sum_{s\in S}\|\one-s\|_{x'}^2=2\cdot x'(\Delta).
\end{equation}
Every non-trivial cocycle can be normalized by scaling it.

To prove that $\Gamma$ has property~(T), one hence has to show that there does not exist a function $\bar{x}'$ conditionally of positive type with $x'(\Delta)=1$ and $x'(\Delta^2)=0$.
The key point is that this can be done by only looking at finite subsets of $\Gamma$.

Let $\mathcal{T}$ denote the collection of finite subsets $T\subset\Gamma$ that contain $\{\one\}\cup S$ and are connected in the Cayley graph of $\Gamma$.
For any $T\in\mathcal{T}\cup\{\Gamma\}$ let $\mathcal{F}'_T$ be the set of functionals $x'\colon \I[T^{-1}T]\to\R$ for which
$x'(\Delta)=1$ and for which $\langle\cdot,\cdot\rangle_{x'}\vcentcolon=x'\circ\Phi$ is a positive semidefinite form on $\I[T]$.
In the same way as for full functionals $\I[\Gamma]\to\R$, we obtain a Hilbert space $\mathcal{H}_{x'}\vcentcolon=\I[T]/\{\|\cdot\|_{x'}=0\}$ and a map $c_{x'}\colon T\to\mathcal{H}_{x'}$, $\gamma\mapsto\one-\gamma$ such that
$\|\gamma\gamma_1-\gamma\gamma_2\|_{x'}=\|\gamma_1-\gamma_2\|_{x'}$ for all $\gamma,\gamma_1,\gamma_2\in\Gamma$ where the expression makes sense. We think of $c_{x'}$ as a \emph{partial cocycle}, encoded by the partial function $\bar{x}'$.

\begin{lemma}\label{lem:cocycle-accumulation}\leavevmode
\begin{enumerate}
\item
For any $T\in\mathcal{T}$, $x'\in\mathcal{F}'_T$ and $\gamma_1,\gamma_2\in T$ the value $|\bar{x}'(\gamma_1^{-1}\gamma_2)|^{1/2}$ is bounded by the distance between $\gamma_1$ and $\gamma_2$ in the partial Cayley graph $(T,S)\subset(\Gamma,S)$.
\item
For any $T\in\mathcal{T}$ the infimum $\inf_{x'\in\mathcal{F}'_T}x'(\Delta^2)$ is assumed.
\item
$\inf_{x'\in\mathcal{F}'_\Gamma}x'(\Delta^2)=\sup_{T\in\mathcal{T}}\min_{x'\in\mathcal{F}'_T}x'(\Delta^2)$
and the infimum is assumed.
\end{enumerate}
\end{lemma}

\begin{proof}
\textbf{1)}
For any $t_1,t_2\in T$ that are neighbors in the Cayley graph, $t_2=t_1s$, we have
$\|t_1-t_2\|^2_{x'}=\|\one-s\|^2_{x'}\leq\sum_S\|\one-s\|^2_{x'}=2$.
Now, if $\gamma_1=t_0,t_1,\dots,t_n=\gamma_2\in T$ form a path in the graph $(T,S)$, connecting $\gamma_1$ and $\gamma_2$, then by the triangle inequality
\[\sqrt{2}\cdot|\bar{x}'(\gamma_1^{-1}\gamma_2)|^{1/2}=\|\one-\gamma_1^{-1}\gamma_2\|_{x'}=\|\gamma_1-\gamma_2\|_{x'}\leq\sum_{i=1}^n\|t_i-t_{i-1}\|_{x'}\leq\sqrt{2}\cdot n.\]

\textbf{2)} It follows from part 1) that the space
$\mathcal{F}'_T\subset\{\text{functionals on }\I[T^{-1}T]\}$ is compact.

\textbf{3)} First, we note that $\inf_{x'\in\mathcal{F}'_{T_1}}\bar{x}'(\Delta^2)\leq \inf_{x'\in\mathcal{F}'_{T_2}}\bar{x}'(\Delta^2)$ when $T_1\subset T_2$, since every $x'\in\mathcal{F}'_{T_2}$ can be restricted to give an element in $\mathcal{F}'_{T_1}$ with the same curvature. Therefore, the left expression in the statement is at least as big as the right, and in the right expression we may restrict ourselves to a cofinal sequence $T_1\subset T_2\subset\dots$, $\bigcup_i T_i=\Gamma$.

Let now $x'_k\in\mathcal{F}'_{T_k}$ such that $x'_k(\Delta^2)=\min_{x'\in\mathcal{F}'_{T_k}}x'(\Delta^2)$.
We extend $\bar{x}'_k$ to all of $\Gamma$ by mapping $\Gamma\setminus T_k^{-1}T_k$ to $0$. Then, by part 1), the sequence $\bar{x}'_k(\gamma)$ is bounded for all $\gamma\in\Gamma$, and hence the sequence $(\bar{x}'_k)_{k\in\N}$ has an accumulation point $\bar{x}'_\infty$ in the topology of pointwise convergence. This accumulation point is a function conditionally of positive type, satisfying $x'_{\infty}(\Delta)=1$ and $x'_{\infty}(\Delta^2)=\lim x'_k(\Delta^2)=\sup_{T\in\mathcal{T}}\min_{x'\in\mathcal{F}'_T}x'(\Delta^2)$.
\end{proof}

\begin{corollary}\label{cor:dual-t}
$\Gamma$ has a non-trivial harmonic cocycle exactly if
$\min_{x'\in\mathcal{F}'_T}x'(\Delta^2)=0$
for all $T\in\mathcal{T}$.
\end{corollary}

\section{SDP duality}
In this section we show that the task of calculating $\min_{x'\in\mathcal{F}'_T}x'(\Delta^2)$ for a fixed set $T\in\mathcal{T}$ can be phrased as a semidefinite program, and that this SDP is dual to the SDP that was introduced in \cite{NT} to search for a property~(T) witness satisfying Equation~\ref{eqn:ozawa-witness}.

\subsection{Preliminaries on SDP duality}
In semidefinite programming -- and more generally in conic programming -- every optimization problem has a dual problem.
In the literature (see e.g. the textbook \cite{BN}*{Section 2.3}) the dual program is often defined by a concrete manipulation of the vectors and matrices that encode the primal SDP. Here, we prefer to give a slightly different but equivalent formulation that is basis-free and highlights the symmetry between the primal and the dual program.

Let $A\colon V\to W$ be a surjective linear map between finite-dimensional vector spaces, let $w_1,w_2\in W$ be vectors and let $\mathcal{K}\subset V$ be a convex cone, the ``cone of positive elements in $V$''. This is exactly the data needed to describe a conic program. We can visualize it in a commutative diagram as follows.

\begin{equation}\label{eqn:sdp-diagram}\begin{tikzcd}
A^{-1}(\Span\{w_1,w_2\})\arrow[r,hookrightarrow,"i_V"]\arrow[d,twoheadrightarrow,"A"] & V\arrow[d,twoheadrightarrow,"A"]&[-1.2cm]\supset\mathcal{K}\\
\Span\{w_1,w_2\}\arrow[r,hookrightarrow,"i_W"] & W&
\end{tikzcd}\end{equation}

The primal problem considers all elements $w\in\Span\{w_1,w_2\}$ that are of the form $w=w_1+\varepsilon w_2$ and that have a lift $\tilde{w}\in \mathcal{F}\vcentcolon=\{v\in A^{-1}(\Span\{w_1,w_2\})\mid i_V(v)\in\mathcal{K}\}$
satisfying $A(\tilde{w})=w$.
Among these elements one wishes to minimize $\varepsilon=\varepsilon_{\tilde{w}}$.
One calls $\mathcal{F}$ the set of \emph{primal feasible points} and $\mathcal{F}\cap i_V^{-1}(\interior\mathcal{K})$ the set of \emph{primal strictly feasible points}, where $\interior\mathcal{K}$ denotes the interior.

To obtain the dual problem one dualizes all spaces and maps, and exchanges the roles of $w_1,w_2$. Namely, one considers all functionals $w'\colon\Span\{w_1,w_2\}\to\R$ that are of the form $w'=\delta_{w_2}+\varepsilon'\delta_{w_1}$ and that have an extension $\tilde{w}'\in\mathcal{F}'\vcentcolon=\{w'\colon W\to\R\mid w'\circ A\in\mathcal{K}'\}$ satisfying $\tilde{w}'\circ i_W=w'$, where $\mathcal{K}'\vcentcolon=\{\ell\colon V\to\R\mid\ell_{|\mathcal{K}}\geq 0\}$ is the dual cone to $\mathcal{K}$.
Among these functionals one wants to minimize $\varepsilon'=\varepsilon'_{\tilde{w}'}$.
One calls $\mathcal{F}'$ the set of \emph{dual feasible points} and $\mathcal{F}'\cap {A'}^{-1}(\interior\mathcal{K}')$ the set of \emph{dual strictly feasible points}.

In this formulation the symmetry between the primal and the dual problem is obvious: By the Hyperplane Separation Lemma $(\mathcal{K}')'\vcentcolon=\{v\in V\mid(\ev_v)_{|\mathcal{K}'}\geq 0\}$ is the same as $\mathcal{K}$, whence the dual of the dual program is again the primal program when one identifies the vector spaces with their double duals via evaluation.
If we choose bases for the vector spaces and write all maps as matrices, we see that via an affine transformation the above definition of the dual program is equivalent to the definition in \cite{BN}.

\medskip

One immediate consequence of the above definition for the dual program is the inequality
$0\leq(\tilde{w}'\circ A)(i_V(\tilde{w}))
=(\tilde{w}'\circ i_W)(A(\tilde{w}))
=w'(w)=\varepsilon_{\tilde{w}}+\varepsilon'_{\tilde{w}'}$.
Letting $\varepsilon_{\textit{inf}}\vcentcolon=\inf_{\tilde{w}\in\mathcal{F}}\varepsilon_{\tilde{w}}$ and $\varepsilon'_{\textit{inf}}\vcentcolon=\inf_{\tilde{w}'\in\mathcal{F}'}\varepsilon'_{\tilde{w}'}$, the number $\varepsilon_{\textit{inf}}+\varepsilon'_{\textit{inf}}$ is called \emph{duality gap}, and the preceding inequality is called \emph{weak duality}. \emph{Strong duality} holds when the duality gap is zero. The fact that strong duality does hold under mild assumptions on the optimization problem is a central result in conic optimization.

\begin{theorem}[Slater's Constraint Qualification, see \cite{BN}*{Theorem 2.4.1}]\label{thm:strong-duality}
Assume that the primal conic program has a non-empty set of strictly feasible points and is bounded below in the sense that $\inf_{\tilde{w}\in\mathcal{F}}\varepsilon_{\tilde{w}}>-\infty$.

Then the dual program attains its infimum $\varepsilon'_{\textit{inf}}=\inf_{\tilde{w}'\in\mathcal{F}'}\varepsilon'_{\tilde{w}'}$ at some optimal $\tilde{w}'$ and the duality gap vanishes.
\end{theorem}

Since the preceding theorem lies at the heart of this article, we want to assure the reader that its proof is accessible without having to wade through the territories of optimization theory.
In fact, with the above formulation of the dual program, it can be viewed as a direct consequence of the M.\ Riesz Extension Theorem.

\begin{theorem}[M.\ Riesz Extension Theorem]\label{thm:extension-theorem}
Let $Y$ be a vector space, $X\subset Y$ a subspace, and $K\subset Y$ a convex cone. We say that a functional $\ell\colon X\to\R$ or $Y\to\R$ is \emph{positive} if $\ell_{|K}\geq 0$.

If $Y=X+K$, then any positive functional on $X$ can be extended to a positive functional on $Y$.
\end{theorem}

\begin{proof}[Proof of Theorem \ref{thm:strong-duality}]
Consider the convex cone $K\vcentcolon=A(\mathcal{K})\subset W$. Its restriction to $\Span\{w_1,w_2\}$ looks like a pie segment in $\R^2$, and the point $w_1+\varepsilon_{\textit{inf}}\cdot w_2$ lies on its boundary. Hence, the functional $\ell\vcentcolon=\delta_{w_2}-\varepsilon_{\textit{inf}}\cdot\delta_{w_1}$ is positive on $K\cap\Span\{w_1,w_2\}$ and $\inf_{\tilde{w}\in\mathcal{F}}\ell(w)=0$. The primal SDP having strictly feasible points means that $\Span\{w_1,w_2\}$ intersects the interior of $K$, whence $\Span\{w_1,w_2\}+K=W$. By Theorem~\ref{thm:extension-theorem} we can extend $\ell$ to a functional $\tilde{w}'$ on $W$ such that $\tilde{w}'\circ A$ is positive on $\mathcal{K}$.
\end{proof}

\subsection{Duality for the property~(T) SDPs}\label{sec:T-duality}
We now specialize the conic program of Diagram~\ref{eqn:sdp-diagram} to the property~(T) SDPs. Using the terminology of Section~\ref{sec:cocycle-perspective}, we set for a fixed $T\in\mathcal{T}$
\begin{align*}
V&\vcentcolon=\I[T]\otimes\I[T],&A&\vcentcolon=\Phi\colon\gamma_1\otimes\gamma_2\mapsto\gamma_1^{-1}\gamma_2,&w_1&\vcentcolon=\Delta^2,\\
W&\vcentcolon=\I[T^{-1}T],&
\mathcal{K}&\vcentcolon=\big\{\textstyle\sum\omega_i\otimes\omega_i\,\big\vert\,\omega_i\in\I[T]\big\},&w_2&\vcentcolon=\Delta.
\end{align*}
Then the primal SDP describes the task to write $\Delta^2+\varepsilon\Delta$ as a sum $\sum_i\omega_i^*\omega_i$ for $\omega_i\in\I[T]$ and $\varepsilon$ minimal. To prove property~(T) one has to show that the optimal $\varepsilon$ is $<0$.
The dual SDP is to compute $\min_{x'\in\mathcal{F}'_T}x'(\Delta^2)$, that is, to find a normalized partial cocycle on $T$ with minimal curvature $\varepsilon'$. To prove property~(T) one has to show that the optimal $\varepsilon'$ is $>0$.

\begin{lemma}\label{lem:prop-t-duality}
Let $T\in\mathcal{T}$. Then $\inf_{x\in\mathcal{F}_T}\delta_\Delta\circ\Phi(x)=-\inf_{x'\in\mathcal{F}'_T}x'(\Delta^2)$ and both infima are assumed.
\end{lemma}
\begin{proof}
The restriction of the trace $\tr\colon\R[\Gamma]\to\R$, $\sum\lambda_i\gamma_i\mapsto\lambda_\one$ to $\I[T^{-1}T]$ is strictly positive and satisfies $\tr(\Delta)\neq 0$. By scaling it one obtains a strictly feasible point for the dual SDP. The dual SDP is also bounded, since $\|\Delta\|_{x'}^2\geq 0$.
Now we apply Theorem~\ref{thm:strong-duality} to the dual SDP, using the fact that the dual of the dual SDP is again the primal SDP.
This shows that the duality gap vanishes and that the infimum $\inf_{x\in\mathcal{F}_T}\delta_\Delta\circ\Phi(x)$ is assumed.
By part 2) of Lemma~\ref{lem:cocycle-accumulation}, the infimum $\inf_{x'\in\mathcal{F}'_T}x'(\Delta^2)$ is also assumed.
\end{proof}

\begin{remark}\label{rmk:primal-strictly-feasible}
The fact that the infimum $\inf_{x'\in\mathcal{F}'_T}x'(\Delta^2)$ is assumed can also be shown by applying Theorem~\ref{thm:strong-duality} to the primal SDP, which has a strictly feasible point because $\Delta$ is an order unit for $\I[T]$ \cite{O}*{Lemma 2}.
\end{remark}

By combining Lemma~\ref{lem:prop-t-duality} with Corollary~\ref{cor:dual-t} one gets a direct proof for Proposition~\ref{prp:dichotomy}. Using Ozawa's theorem \cite{O}*{Prop.\ 4}, we also get a quantitative result.

\begin{theorem}\label{thm:quantitative-version}
Let $\Gamma$ be a finitely generated group, $S=S^{-1}$ a finite generating set and $\Delta\in\R[\Gamma]$ the corresponding Laplacian.
Let $\mathcal{C}$ denote the set of 1-cocycles $\Gamma\to\mathcal{H}$ satisfying $\sum_{s\in S}\|c(s)\|^2=2$ for any orthogonal representation $\Gamma\curvearrowright\mathcal{H}$.

Then the infimum $\inf_{c\in\mathcal{C}}\|\sum_{s\in S}c(s)\|^2$ is assumed and agrees with the size of the spectral gap of $\Delta$ above $0$.
\end{theorem}
\begin{proof}
Using the correspondence of cocycles with $\Gamma$-invariant positive semidefinite forms on $\I[\Gamma]$, we may replace $\mathcal{C}$ with $\{c_{x'}\mid x'\in\mathcal{F}'_\Gamma\}$, and we recall that
$\|\sum_{s\in S}c_{x'}(s)\|^2=\|\Delta\|_{x'}^2=x'(\Delta^2)$
for all $x'\in\mathcal{F}'_\Gamma$.
By Lemma~\ref{lem:cocycle-accumulation}, $\inf_{x'\in\mathcal{F}'_\Gamma}x'(\Delta^2)$ is assumed and equals
$\sup_{T\in\mathcal{T}}\inf_{x'\in\mathcal{F}'_T}x'(\Delta^2)$. By Lemma~\ref{lem:prop-t-duality}, this is the same as $-\sup_{T\in\mathcal{T}}\inf_{x\in\mathcal{F}_T}\delta_\Delta\circ\Phi(x)$, and this value agrees with the size of the spectral gap of $\Delta$ by \cite{O}*{Prop.\ 4}.
\end{proof}

We emphasize that in practice the cocycle SDP leads to exactly the same computer calculation as the SDP for the spectral gap witness search.
To prove that $\inf_{x'\in\mathcal{F}'_T}x'(\Delta^2)$ is $>0$, one shows that $\inf_{x\in\mathcal{F}_T}\delta_\Delta\circ\Phi(x)<0$ and applies weak duality.
However, the fact that strong duality holds makes it possible to go back and forth between the algebraic primal formulation and the geometric dual formulation of the SDP.

\section{Going back and forth between the primal and the dual property~(T) SDP}\label{sec:dualizing-arguments}
We now illustrate how SDP duality can be used to translate between the geometric dual perspective and the algebraic primal perspective of the property~(T) SDP, and how this can be used to better understand and to simplify the SDP.

\subsection{Bounds from the triangle and parallelogram inequalities}\label{sec:triangle-parallelogram}
First, we consider the triangle inequality and the parallelogram inequality, which are the dual of the algebraic arguments in \cite{O}*{Lemma 2}, respectively \cite{KKN}*{Lemma 3.6}.

To prove that $\Gamma$ has property~(T), one wants to bound the distances $\|c_{x'}(\gamma)\|=\|\one-\gamma\|_{x'}$ from above, uniformly for all $x'\in\mathcal{F}'$.
Sufficiently good bounds in the $2$-ball $\gamma\in S^2\subset\Gamma$ would imply that the mean $1/|S|\cdot\sum_S c_{x'}(s)$ must always have strictly positive distance from $0=c_{x'}(\one)$, whence no harmonic cocycles exist.
The normalizing condition from Equation~\ref{eqn:normalizing-condition} implies $\|\one-s\|_{x'}^2\leq 2$ for $s\in S$. In some cases, like $\SL(n,\Z)$ \cite{KNO}, one may additionally assume that $x'$ carries a symmetry that implies $\|\one-s_1\|_{x'}=\|\one-s_2\|_{x'}$ for all $s_1,s_2\in S$, whence $\|\one-s\|=\sqrt{2/|S|}$.

The easiest way to obtain further bounds for any group $\Gamma$ is by the triangle inequality. If one already has bounds $\|\one-\gamma_i\|_{x'}\leq\mu_i$ for $\gamma_1,\gamma_2\in T$, $x'\in\mathcal{F}'_T$, then
\[\|\one-\gamma_1\gamma_2\|_{x'}\leq\|\one-\gamma_1\|_{x'}+\|\gamma_1-\gamma_1\gamma_2\|_{x'}\leq\mu_1+\mu_2.\]
To transfer this argument to the primal SDP, we observe that
\begin{align*}
\mathrm{Trg}_{\gamma_1,\gamma_2}\colon x'&\mapsto-\|\one-\gamma_1\gamma_2\|_{x'}^2+(\mu_1+\mu_2)^2\cdot x'(\Delta)\\
&\phantom{\vcentcolon}=-x'\big((\one-\gamma_1\gamma_2)^*(\one-\gamma_1\gamma_2)\big)+(\mu_1+\mu_2)^2\cdot x'(\Delta)
\end{align*}
is a linear functional on the dual space of $\I[T^{-1}T]$, and positive with respect to the dual cone of $\Phi(\mathcal{K})$. By the Hyperplane Separation Lemma, it then follows that the element $\omega\vcentcolon=-(\one-\gamma_1\gamma_2)^*(\one-\gamma_1\gamma_2)+(\mu_1+\mu_2)^2\cdot \Delta$ lies in $\Phi(\mathcal{K})$ and can hence be used as part of the sum in the decomposition of Equation~\ref{eqn:ozawa-witness}.

In this simple case it is easy to construct a decomposition $\omega=\sum x_i^*x_i$ that directly shows that $\omega$ is positive.
The edge-case where the inequalities \[\|\one-\gamma_1\gamma_2\|_{x'}\leq\|\one-\gamma_1\|_{x'}+\|\one-\gamma_2\|_{x'}\leq\mu_1+\mu_2\]
become equalities occurs when the bounds $\mu_i$ are sharp and the points $c_{x'}(\one)$, $c_{x'}(\gamma_1)$, $c_{x'}(\gamma_2)$ lie on a line, with $\mu_2 c_{x'}(\gamma_1)+\mu_1 c_{x'}(\gamma_2)=0$. That happens exactly when the element $x\vcentcolon=\mu_2\cdot(\one-\gamma_1)+\mu_1\cdot(\one-\gamma_2)$ satisfies $\|x\|_{x'}^2=0$. This leads us to consider
\begin{align*}
x^*x&=({\mu_2}^2+\mu_1\mu_2)\cdot(\one-\gamma_1)^*(\one-\gamma_1)+({\mu_1}^2+\mu_1\mu_2)\cdot(\one-\gamma_2)^*(\one-\gamma_2)\\
&\phantom{=}-\mu_1\mu_2\cdot(\one-\gamma_1\gamma_2)^*(\one-\gamma_1\gamma_2).
\end{align*}
By assumption, $\omega_1\vcentcolon={\mu_1}^2\Delta-(\one-\gamma_1)^*(\one-\gamma_1)$ and $\omega_2\vcentcolon={\mu_2}^2\Delta-(\one-\gamma_2)^*(\one-\gamma_2)$ are both positive, so $\omega=(\mu_1\mu_2)^{-1}\cdot x^*x+(\mu_2\mu_1^{-1}+1)\cdot\omega_1+(\mu_1\mu_2^{-1}+1)\cdot\omega_2$ is also positive.

\begin{remark}
The elements $x^*x$ for $\mu_1=\mu_2$ lie at the core of the proof that $\Delta$ is an order unit for $\I[\Gamma]$ \cite{O}*{Lemma 2} and for the error estimate \cite{NT}*{Lemma 2.1}. The statements that $\Delta$ is an order unit and that each distance $\|\one-\gamma\|_{x'}$ is uniformly bounded over all $x'$ with $x'(\Delta)=1$ are dual to each other.
The error estimates of \cite{NT}*{Lemma 2.1} can be slightly improved by using the above decomposition and allowing $\mu_1\neq\mu_2$.
\end{remark}

The next more difficult bounds beyond the triangle inequality are obtained from the geometric relations between four points $0,\one-\gamma_1,\one-\gamma_2,\one-\gamma_3$. For simplicity, we restrict to the case where we know that $\|\one-\gamma_1\|_{x'}=\|\gamma_2-\gamma_3\|_{x'}$ and $\|\one-\gamma_2\|_{x'}=\|\gamma_1-\gamma_3\|_{x'}$. One situation where this happens is when $\gamma_1$ and $\gamma_2$ commute and $\gamma_3=\gamma_1\gamma_2$. Then, the parallelogram inequality says that $\|\one-\gamma_3\|_{x'}^2+\|\gamma_1-\gamma_2\|_{x'}^2\leq 2\cdot\|\one-\gamma_1\|^2+2\cdot\|\one-\gamma_2\|_{x'}^2$. The inequality becomes an equality when all four points lie on a plane, satisfying the dependence relation $\one-\gamma_3=\one-\gamma_1+\one-\gamma_2$. This happens exactly when $x\vcentcolon=\one-\gamma_1-\gamma_2+\gamma_3$ satisfies $\|x\|_{x'}^2=0$. As before, one can use the summand $x^*x$ to obtain the analog of the parallelogram inequality on the primal side. The element $x^*x$ plays a crucial role in the proof of \cite{KKN}*{Lemma 3.6}.

\subsection{A condition on harmonic cocycles}\label{sec:main-simplification}
Next, we turn to the simplification that will be the basis for our computation on $\Aut(F_4)$ below. Here, coming from the opposite direction as before, we find a condition that all harmonic cocycles must necessarily satisfy. We then use this condition as a restriction in the SDP search for a harmonic cocycle.

Let $(\Gamma,S)$ be as before and let $H'<\Aut(\Gamma)$ be the subgroup of automorphisms that leave $S$ invariant. Since $S$ is a finite generating set, $H'$ is finite. By averaging its associated function conditionally of positive type over $H'$, one can make every harmonic cocycle $H'$\nobreakdash-invariant, meaning that $\|\one-\gamma\|_{x'}=\|\one-h'(\gamma)\|_{x'}$ for all $\gamma\in\Gamma$, $h'\in H'$. The analog of this argument on the primal side was one key insight of \cite{KNO}.

For inner automorphisms one can say even more. Let $H<\Gamma$ be the preimage of $H'$ under the conjugation homomorphism $\mathrm{conj}\colon\Gamma\to\mathrm{Inn}(\Gamma)<\Aut(\Gamma)$. Note that $H$ is finite if $\Gamma$ has trivial center.

\begin{lemma}\label{lem:main-simplification}
Let $\Gamma,S,H$ be as above and let $c_{x'}$ be a harmonic cocycle for $(\Gamma,S)$.
Then the right multiplication action of any $h\in H$ on $\I[\Gamma]$ induces on $\I[\Gamma]/\{\|\cdot\|_{x'}=0\}$ a translation by $h-\one$. If $h$ has finite order, then $\|\one-\gamma h\|_{x'}=\|\one-\gamma\|_{x'}$ for all $\gamma\in\Gamma$, in particular $\|\one-h\|_{x'}=0$.
\end{lemma}
\begin{proof}
Let $\langle\cdot,\cdot\rangle=\langle\cdot,\cdot\rangle_{x'}$.
The idea is that the two sets $\{\one-s\}_{s\in S},\{\one-sh\}_{s\in S}\subset\I[\Gamma]$ satisfy $\|(\one-s)-(\one-sh)\|=\|\one-h\|\,\forall s\in S$, and hence the distance between their means achieves its maximum $\|\one-h\|$ only when the sets are translates of each other modulo $\{\|\cdot\|=0\}$.

Formally, it follows from $\|\Delta\|^2=0$ that $|\langle\xi,\eta\Delta\rangle|\leq\|\eta^{-1}\xi\|\cdot\|\Delta\|=0$ for any $\xi\in\I[\Gamma]$, $\eta\in\Gamma$. Applying this to $\xi=\gamma-\gamma h$, we get
\begin{align*}
|S|\cdot\langle\gamma-\gamma h,\gamma-\gamma h\rangle
&=\langle\gamma-\gamma h,|S|\cdot\gamma-\gamma\Delta-|S|\cdot\gamma h+\gamma h\Delta\rangle\\
&=\big\langle\gamma-\gamma h,\big({\textstyle\sum_{s\in S}}\gamma s\big)-\big({\textstyle\sum_{s\in S}}\gamma hs\big)\big\rangle\\
&=\big\langle\gamma-\gamma h,\big({\textstyle\sum_{s\in S}}\gamma s\big)-\big({\textstyle\sum_{s\in S}}\gamma sh\big)\big\rangle\\
&=\sum_{s\in S}\langle\gamma-\gamma h,\gamma s-\gamma sh\rangle\\
&\leq\sum_{s\in S}\|\gamma-\gamma h\|\cdot\|\gamma s-\gamma sh\|\\
&=\sum_{s\in S}\|\gamma-\gamma h\|\cdot\|\gamma-\gamma h\|\\
&=|S|\cdot\langle\gamma-\gamma h,\gamma-\gamma h\rangle.
\end{align*}
Since the first term equals the last, the inequality must be an equality. This can only happen when for every $s\in S$ the vector $\gamma s-\gamma sh$ is parallel to $\gamma-\gamma h$ in $\I[\Gamma]/\{\|\cdot\|=0\}$. Since $\gamma$ was arbitrary and $S$ generates $\Gamma$, it follows that $\gamma-\gamma h$ is parallel to $\one-h$ in $\I[\Gamma]/\{\|\cdot\|=0\}$ for all $\gamma\in\Gamma$.

For the second statement, if $h^m=\one$, then it follows from the first part that
\[0=\|{\textstyle\sum_{i=1}^m}h^{i-1}-h^i\|=\sum_{i=1}^m\|h^{i-1}-h^i\|=m\cdot\|\one-h\|,\]
and, by the reverse triangle inequality,
\[\big|\|\one-\gamma h\|-\|\one-\gamma\|\big|\leq\|\one-\gamma h-(\one-\gamma)\|=\|\gamma-\gamma h\|=\|\one-h\|=0.\qedhere\]
\end{proof}

\begin{example}
Let $\Gamma=H^3\vcentcolon=\langle e,f,g\mid [e,f]=[f,g]=\one,[e,g]=f\rangle$ be the Heisenberg group, $S\subset H^3$ any finite symmetric generating set, and $c_{x'}$ a harmonic cocycle for $(\Gamma,S)$. The left multiplication action of $f$ on $\I[\Gamma]$ coincides with its right action, which by Lemma~\ref{lem:main-simplification} induces a translation on $\I[\Gamma]/\{\|\cdot\|_{x'}=0\}$. Therefore, $\|\one-f^{k^2}\|_{x'}=\|\sum_{i=1}^{k^2} f^{i-1}-f^i\|_{x'}=k^2\cdot\|\one-f\|_{x'}$.
But by the triangle inequality we have $\|\one-f^{k^2}\|_{x'}=\|\one-e^k g^k e^{-k} g^{-k}\|_{x'}\leq 2k\cdot\|\one-e\|_{x'}+2k\cdot\|\one-g\|_{x'}$ and, letting $k\to\infty$, we conclude that $f$ must act trivially on $\mathcal{H}_{x'}$.

This means that any harmonic cocycle on $H^3$ factors through the quotient map $q\colon H^3\to H^3/\langle f\rangle\cong\Z^2$ postcomposed with a harmonic cocycle on $(\Z,q(S))$. As $\Z^2$ is abelian, that cocycle must come from an action $\Z^2\curvearrowright\mathcal{H}_{x'}$ by translations.
\end{example}

We can use Lemma~\ref{lem:main-simplification} to simplify the property~(T) SDP.
In the dual formulation, we now search for a cocycle $x'\in\mathcal{F}'$ with minimal curvature $\|\Delta\|_{x'}^2$ that also satisfies the additional constraints $\|\one-\gamma\|_{x'}^2=\|\one-\gamma h\|_{x'}^2$ for all $\gamma\in\Gamma$ and $h\in H$ of finite order. Concretely, this means to divide the space $W_0\vcentcolon=\Span\{\gamma-\gamma h\mid \gamma\in\Gamma, h\in H \text{ of finite order}\}$ out of $W=\I[T^*T]$ in Diagram~\ref{eqn:sdp-diagram}. It also means that we need only one representative from each $H$-coset in the supporting set $T$, since additional representatives do not contribute more information.

In the primal formulation of the SDP, the effect of the simplification is that we search for a decomposition $\Delta^2-\varepsilon\Delta\equiv\sum_{i=1}^n \omega_i^*\omega_i\mod W_0$.

\begin{remark}\label{rem:primal-simplification}
By dualizing every step in the proof of Lemma~\ref{lem:main-simplification} in a similar way as in Section~\ref{sec:triangle-parallelogram}, one obtains, as the corresponding result on the side of the primal SDP, an explicit procedure to turn a decomposition $\Delta^2-\varepsilon\Delta\equiv\sum_{i=1}^n \omega_i^*\omega_i\mod W_0$, $\varepsilon>0$, into a decomposition $N\Delta^2-\varepsilon'\Delta=\sum_{i=1}^n \hat\omega_i^*\hat\omega_i$, $N\gg 0$, $\varepsilon'>0$. This procedure involves many auxiliary $\varepsilon$ and we do not include it here.
\end{remark}

\section{Application to \texorpdfstring{$\SL(n,\mathbb{Z})$}{SL(n,Z)}}

We now discuss the example of $\Gamma=\SL(n,\Z)$, $n\geq 3$, arguably the most prominent example of discrete property~(T) groups, from the perspective of the dual SDP. As the generating set we take $S\vcentcolon=\{e_{ij}^{\pm 1}\mid 1\leq i\neq j\leq n\}$, where $e_{ij}$ is the elementary matrix that differs from $\one$ only by an entry of $1$ in position $(i,j)$. Let $\Delta$ be the associated Laplacian.
We let further $H_0<\GL(n,\Z)$ be the subgroup generated by the permutation matrices and the diagonal matrices, $H\vcentcolon=H_0\cap\SL(n,\Z)$, and we let $H'<\Aut(\SL(n,\Z))$ be the subgroup generated by the inner automorphisms $\{\mathrm{conj}_{h}\mid h\in H_0\}$ and the ``exceptional automorphism'' $X\mapsto (X^\transp)^{-1}$.

The property~(T) dual SDP searches for a $\Gamma$-invariant seminorm $\|\cdot\|\vcentcolon=\|\cdot\|_{x'}$ on $\I[\Gamma]$ that minimizes $\|\Delta\|^2$ under the constraint $\sum_{s\in S}\|\one-s\|^2=2$. By averaging over $H'$, we may assume the $H'$-symmetry $\|\one-h'(\gamma)\|=\|\one-\gamma\|$. In particular, after rescaling, we may assume that $\|\one-s\|=1$ for all $s\in S$. Then $\|\Delta\|^2$ can be rewritten as
\begin{align*}
2\cdot\|\Delta\|^2&=2\cdot\big\langle\,{\textstyle\sum_{s\in S}} \one-s,\,{\textstyle\sum_{s'\in S}}\one-s'\,\big\rangle
=\sum_{s,s'\in S}2\cdot\langle\,\one-s,\one-s'\,\rangle\\
&=\sum_{s,s'\in S} \|\one-s\|^2+\|\one-s'\|^2-\|s-s'\|^2
=2|S|^2-\sum_{s,s'}\|\one-s^{-1}s'\|^2\\
&=2|S|^2
-|S|\cdot\|\one-e_{12}e_{12}\|^2
-|S|\cdot\|\one-e_{12}e_{21}\|^2
-|S|\cdot\|\one-e_{12}e_{21}^{-1}\|^2\\
&\phantom{=}-|S|\cdot4(n-2)\cdot\|\one-e_{12}e_{13}\|^2
-|S|\cdot4(n-2)\cdot\|\one-e_{12}e_{23}\|^2\\
&\phantom{=}-|S|\cdot2(n-2)(n-3)\cdot\|\one-e_{12}e_{34}\|^2,
\end{align*}
using $H'$-symmetry for the last equality.

To show that $\|\Delta\|^2=0$ cannot be achieved, one has to bound the summands on the right-hand side of the equation.
By the parallelogram inequality, $\|\one-e_{12}e_{34}\|^2\leq 2$. Also,
\[
8-\|\one-e_{12}e_{12}\|^2-\|\one-e_{12}e_{21}\|^2-\|\one-e_{12}e_{21}^{-1}\|^2=
1/4\cdot\|4\cdot\one-e_{12}-e_{12}^{-1}-e_{21}-e_{21}^{-1}\|^2\geq 0.
\]
Hence, showing $\|\one-e_{12}e_{13}\|^2+\|\one-e_{12}e_{23}\|^2<4$ would imply property~(T), and even more, it would imply property~(T) also for $\SL(n',\Z)$, $n'\geq n$.

\begin{remark}
This is exactly the dual argument of the strategy of \cite{KKN}. Another way of phrasing it is that, by $H'$-symmetry, the inner products $\langle\Delta_{ij},\Delta_{i'j'}\rangle$ of the group ring elements $\Delta_{ij}\vcentcolon=4\cdot\one-e_{ij}-e_{ij}^{-1}-e_{ji}-e_{ji}^{-1}$ take only three different values corresponding to whether $(i,j)=(i',j')$, or $i,i',j,j'$ are pairwise distinct, or all other cases. The first two values are non-negative and if the last value can be shown to be strictly positive, then $\Delta=\sum_{i<j}\Delta_{ij}$ must have non-trivial seminorm.
\end{remark}

\medskip

The simplification from Section~\ref{sec:main-simplification} now says that $\|\one-\gamma h\|=\|\one-\gamma\|$ can be assumed for all $h\in H$, $\gamma\in\Gamma$. The most visible consequence of this is that $e_{21}^{-1}e_{12}e_{21}^{-1}\in H$, whence $\|\one-e_{12}e_{21}^{-1}\|=\|e_{21}^{-1}-e_{21}^{-1}e_{12}e_{21}^{-1}\|=\|e_{21}^{-1}-\one\|=1$.

\begin{remark}
Consider the Steinberg groups
\[\St(n,\Z)\vcentcolon=\big\langle\, \{E_{ij}\}_{1\leq i\neq j\leq n}\mid [E_{ij},E_{jk}]=E_{ik}\,\forall i\neq k; [E_{ij},E_{kl}]=\one\,\forall i\neq l,j\neq k\,\big\rangle\]
with the generating set $S'\vcentcolon=\big(\St(n,\Z)\to\SL(n,\Z)\big)^{-1}(S)$.
Since conjugation with the non-trivial element in $\ker\big(\St(n,\Z)\to\SL(n,\Z)\big)$ leaves $S'$ invariant \cite{M}*{\S 10}, the simplified property~(T) SDP for $(\St(n,\Z),S')$ is exactly the same as for $(\SL(n,\Z),S)$.
\end{remark}

\begin{figure}[ht]
\centering
\includegraphics[scale=0.75]{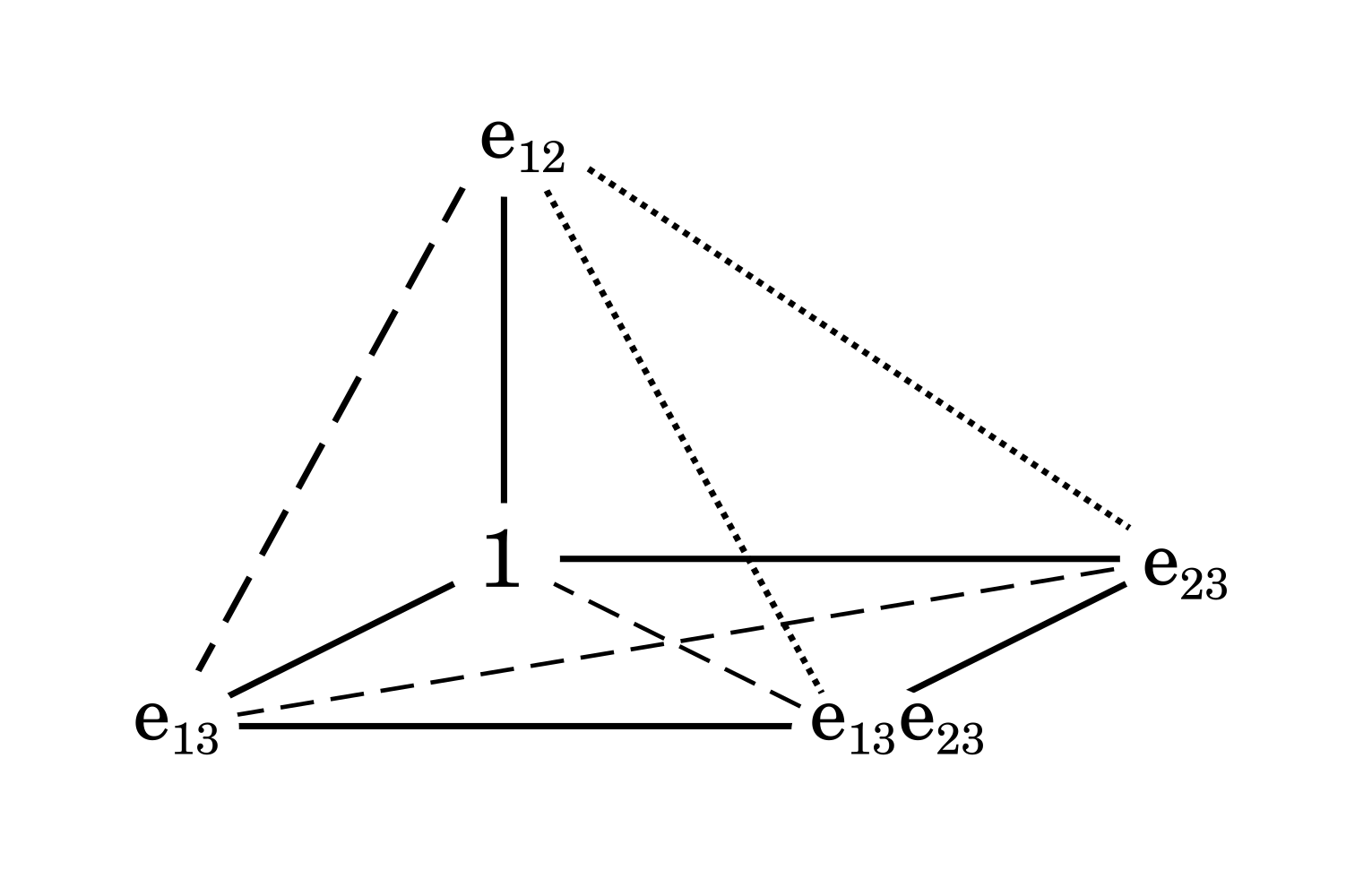}\hspace{2cm}
\includegraphics[scale=0.75]{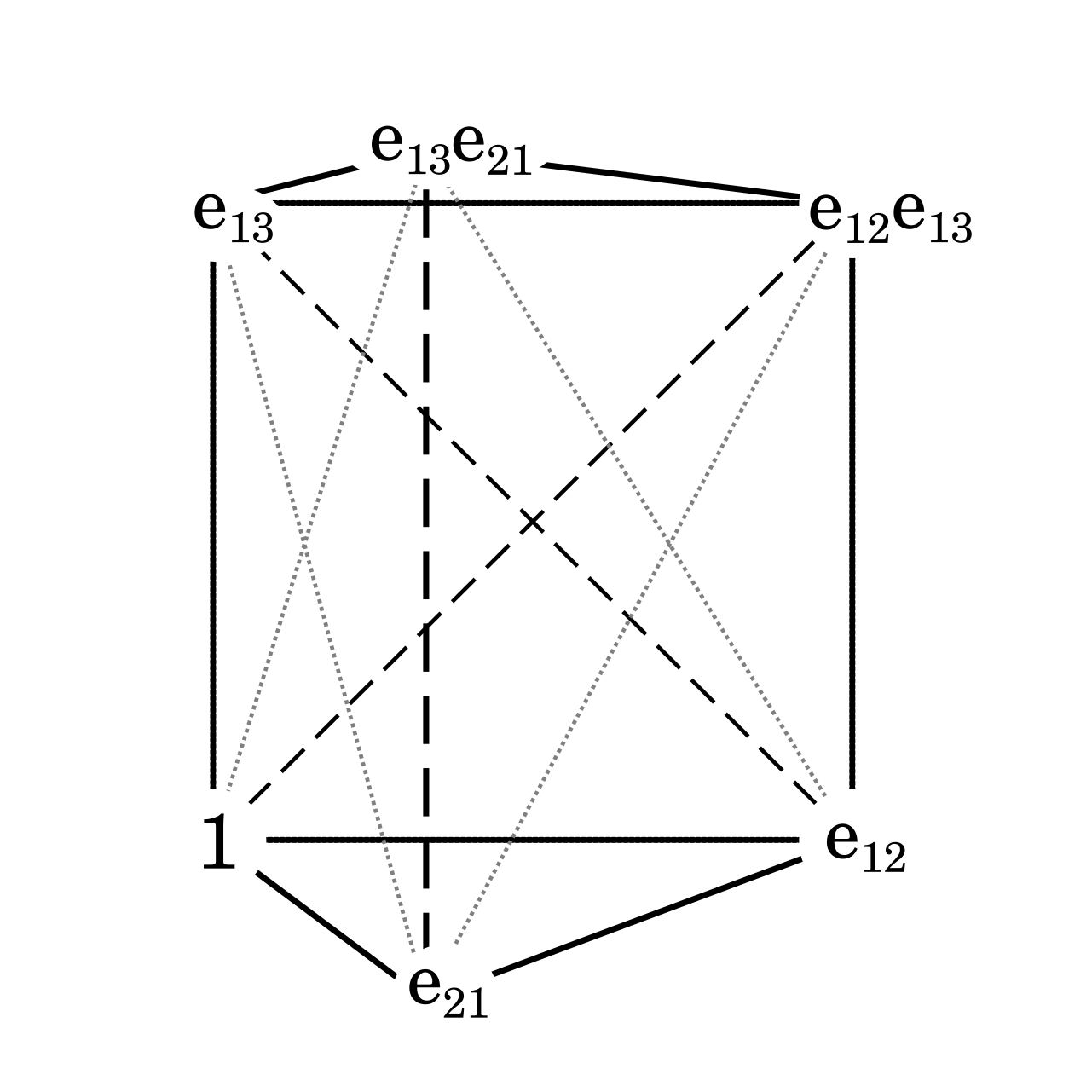}
\caption[Two manually computable bounds]{
In the simplified SDP the distances between the group elements shown in the two diagrams take only three different values, which are drawn as solid, dashed and dotted lines, respectively.}
\label{fig:simple-bounds}
\end{figure}

We attempted to prove $\|\one-e_{12}e_{13}\|^2+\|\one-e_{12}e_{23}\|^2<4$ in the setting of the simplified SDP for $\SL(3,\Z)$.
Easy geometric calculations suffice to prove bounds for $\|\one-e_{12}e_{13}\|^2$ and $\|\one-e_{12}e_{23}\|^2$ that are quite close to being sufficient.
The left diagram in Figure~\ref{fig:simple-bounds} shows that $\|\one-e_{12}e_{13}\|^2$ is strictly less than $2$. It leads to the bound
\begin{equation}\label{eqn:K-bound}
\|\one-e_{12}e_{13}\|^2\leq\text{largest root of }\,2x^3-2x^2-4x+1\approx 1.91.
\end{equation}
The right diagram in Figure~\ref{fig:simple-bounds} leads, in another straightforward calculation, to the bound
\begin{equation}\label{eqn:N-bound}
\|\one-e_{12}e_{23}\|^2\leq 1+\sqrt{2}/2\cdot\|\one-e_{12}e_{13}\|^2.
\end{equation}

We did not manage to prove $\|\one-e_{12}e_{13}\|^2+\|\one-e_{12}e_{23}\|^2<4$ with a completely manual calculation. But we were able to decrease the problem complexity needed for a successful computer proof to a point that is almost in reach of human calculation.
The smallest successful supporting set was a set of $19$ group elements in $\Gamma$, together with the inequality $\|e_{12}+e_{12}^{-1}+e_{23}+e_{23}^{-1}+e_{31}+e_{31}^{-1}-e_{21}-e_{21}^{-1}-e_{32}-e_{32}^{-1}-e_{13}-e_{13}^{-1}\|^2\geq 0$.
To find this set we looked at an optimal dual solution for a larger support $T\subset\Gamma$ and searched among all subsets of $T$ where the maximal distance between any two points was bounded by a chosen value.

We also found other successful supports that were slightly bigger but exhibited high symmetry, e.g.\ by the dihedral group of order $24$.
Using SDPLR \cite{SDPLR}, we found a single vector $\xi\in\I[\Gamma]$, such that $\|\xi\|^2\geq 0$ together with $H'$-symmetry directly implies that $\|\one-e_{12}e_{13}\|^2+\|\one-e_{12}e_{23}\|^2<4$. There did not appear to be a pattern in the coefficients of $\xi$.
Some further empirical observations can be found in an early preprint of this article.\footnote{\url{https://arxiv.org/abs/2009.05134}}

\subsection{Supports in the Heisenberg subgroup and Harper's operator}\label{sec:harpers-operator}
We also observed that the computer can prove $\|\one-e_{12}e_{13}\|^2+\|\one-e_{12}e_{23}\|^2<4$ by using a supporting set $T$ that is contained in the Heisenberg group $H^3<\SL(3,\Z)$. For example,
$T=\big\{{e_{12}}^i {e_{13}}^j {e_{23}}^k\mid i\in\{0,1,2,3,4\}, j\in\{0,1,2,3\}, k\in\{0,1\}\big\}$
suffices in combination with the $H'$-symmetry (even without the additional simplification).
This observation appears meaningful, since it reflects how $\St(n,\Z)$ can be thought of as a number of Heisenberg groups glued together in a symmetric way, and is also reminiscent of how property~(T) can be proved for $\SL(n,\Z)$ via relative property~(T) for $(\SL(2,\Z)\ltimes\Z^2,\Z^2)$.

Motivated by the above observation and the estimate from Equation~\ref{eqn:N-bound}, we studied upper bounds for $\|\one-e_{12}e_{13}\|^2$ that come from looking at the Heisenberg group. Taking the primal perspective, this led us to search for $\eta_1,\eta_2\in\R$, $\eta_1\geq 1$, such that
\begin{align*}
Z&\vcentcolon=-8\cdot\one+(e_{12}+e_{12}^{-1}+e_{23}+e_{23}^{-1})(e_{13}+e_{13}^{-1})\\
&\phantom{\vcentcolon=}+2\eta_1(4\cdot\one-e_{12}-e_{12}^{-1}-e_{23}-e_{23}^{-1})+4\eta_2(2\cdot\one-f-f^*)\in\mathrm{C}^*(H_3)
\end{align*}
is positive and $\eta_1+\eta_2$ is minimal. The element $Z$ is positive, iff its images are positive under the family of $\mathrm{C}^*$-homomorphisms
$\pi_\theta\colon\mathrm{C}^*(H_3)\to\mathcal{B}(L^2(\Z))$, $\theta\in[0,1]$, that map $e_{12}$ to the shift operator, $e_{13}$ to $\exp(2\pi i\theta)\cdot\one$, and $e_{23}$ to multiplication with the function $\exp(2\pi i\theta\,\cdot)$.
Writing this out, the condition that $Z$ is positive is equivalent to
\[\forall\theta\colon H_\theta\vcentcolon= \mathit{Shift}+\mathit{Shift}^*+2\cos(2\pi\theta\,\cdot)\leq\frac{4(\mu_1+\mu_2-1)-4\mu_2\cos\theta}{\mu_1-\cos\theta}\cdot\one.\]

The operator $H_\theta$ is called Harper's operator, a special case of an almost Mathieu operator, which has been much studied in physics. The Hofstadter butterfly arises from plotting its spectrum for rational values of $\theta$. In \cite{BVZ} a norm estimate of $H_\theta$ is used to establish the spectrum of the random walk operator for $H^3$ with generating set $\{{e_{12}}^{\pm 1},{e_{13}}^{\pm 1},{e_{23}}^{\pm 1}\}$.
In \cite{BZ} it is proved that $\|H_\theta\|\leq 2\sqrt{2}$ for $\theta\in[\frac{1}{4},\frac{1}{2}]$ and that
\begin{equation}\label{eqn:Harper-bound}
\|H_\theta\|\leq\sqrt{8+8(\cos\pi\theta-\sin\pi\theta)\cos\pi\theta}\qquad\text{for }\theta\in[0,{\textstyle\frac{1}{4}}].
\end{equation}

In our calculations, this estimate was not good enough to help with the property~(T) proof. Instead, we found that even simple SDPs suffice to locally improve the bound. For example, one can take the supporting set from the left diagram in Figure~\ref{fig:simple-bounds} and consider the positive group ring element
\[0\leq(3\cdot\one-e_{12}-2\cdot e_{13}-2\cdot e_{23}+2\cdot e_{13}e_{23})^*(3\cdot\one-e_{12}-2\cdot e_{13}-2\cdot e_{23}+2\cdot e_{13}e_{23})\in\R[H^3]\]
(the optimal coefficients are slightly more complicated). By averaging this equation over those automorphisms of $H^3$ that map $\{{e_{12}}^{\pm 1},{e_{23}}^{\pm 1}\}$ into itself, and applying the homomorphisms $\pi_\theta$, one obtains the following bound, which is an improvement upon Equation~\ref{eqn:Harper-bound} on the interval from $\theta\approx 0.025$ to $\theta\approx 0.119$.

\begin{proposition}
Harper's operator is bounded by
$\|H_\theta\|\leq\frac{44-40\cdot\cos 2\pi\theta}{13-12\cdot\cos 2\pi\theta}$.
\end{proposition}

Of course, one could obtain even better bounds by solving a suitable SDP on a larger supporting set with the computer.

\begin{remark}
Since the first preprint of this article, Ozawa \cite{O22} was successful in applying the approach of \cite{BZ} to the area of SDP property~(T) proofs, showing that the group $\mathrm{EL}_n(\mathcal{R})$ has a property similar to property~(T) for $n\to\infty$.
\end{remark}

\section{Application to \texorpdfstring{$\Aut(F_4)$}{Aut(F4)}}
The simplification from Section~\ref{sec:main-simplification} leads to a less complex property~(T) SDP and hence faster solving times also for the automorphism groups of the free groups, $\Aut(F_n)$. We used this to prove property~(T) for $\Aut(F_4)$.

Our setup is very similar to \cite{KNO}. Let $F_4=\langle Z\rangle$, $|Z|=4$. We consider the index-$2$ subgroup $\Gamma\vcentcolon=\SAut(F_4)<\Aut(F_4)$, which is the preimage of $\SL(4,\mathbb{Z})$ under the canonical homomorphism $\Aut(F_4)\to\GL(4,\mathbb{Z})$. Since $\SAut(F_4)<\Aut(F_4)$ is a finite index subgroup, it suffices to prove property~(T) for $\SAut(F_4)$. As the generating set for $\SAut(F_4)$ we take the Nielsen generators $S\vcentcolon=\{N_{z_1 z_2}\mid z_1,z_2\in Z\cup Z^{-1}\}$ \cite{G}. Here, $N_{z_1 z_2}\colon F_4\to F_4$ denotes the automorphism that sends $z_1$ to $z_1 z_2$ and is the identity on all other generators in $Z\cup Z^{-1}\setminus\{z_1,z_1^{-1}\}$. Let $\Delta$ be the resulting Laplacian.
We let further $H_0<\Aut(F_4)$ denote the subgroup of those automorphisms that leave $Z\cup Z^{-1}$ invariant. It is generated by the permutations of $Z$ and by the sign flips $z\mapsto z^{\sigma(z)}$ for $\sigma\colon Z\to\{\pm 1\}$. The associated group of inner automorphisms $H'\vcentcolon=\{\mathrm{conj}_{h}\mid h\in H_0\}<\Aut(\SAut(F_4))$ fixes $\Delta$ and can hence be divided out of the property~(T) SDP, as has been done in \cite{KNO}.

\medskip

The reason why the property~(T) proofs of \cites{KNO,KKN} do not cover $\Aut(F_4)$ is that the authors used the $2$-ball $S^2\subset\Gamma$ as the supporting set $T$ in their calculations and -- in contrast to the higher $\Aut(F_n)$ -- this support is apparently not big enough to witness property~(T) for $n=4$. At the same time, the obvious choice for a bigger support, $T=S^3$, leads to an SDP that is much too large to solve in practice.

Our approach was to search for a relatively small but sufficient supporting set by adding to $S^2$ various $H'$-orbits of group elements in $S^3\setminus S^2$. To do so, we needed a heuristic for how much closer the addition of a single orbit brought us to a successful property~(T) proof. The output of the usual property~(T) SDP does not help with this, since it reports a witnessed spectral gap size of exactly $0$ whenever the proof is unsuccessful. Instead, we ran an SDP that was searching for a spectral gap witness for the Laplacian of a higher $\Aut(F_n)$, but using only the given supporting set $T\subset\Aut(F_4)$ together with $H'$-symmetry. This is almost the same SDP that Kaluba--Kielak--Nowak used in \cite{KKN}, concretely it is obtained by introducing suitable coefficients into the decomposition of $\Delta^2$ in \cite{KKN}*{Lemma 3.5}.
The modified SDP produced non-zero spectral gap estimates that we used as a heuristic. Finally, we built our supporting set from a combination of those $H'$-orbits that appeared the most helpful in relation to how much their addition increased the problem size.

\medskip

In practice, the search for a good supporting set depended on a reasonably fast implementation of the property~(T) SDP. We divided out the $H'$-symmetry of both the primal semidefinite variable and the dual variables, as was first done in \cite{KNO}. To solve the SDP we followed \cite{NT} and used the SeDuMi solver \cite{sedumi}, which works with the interior-point method. To construct the SDP, we implemented it directly instead of using a library.

Finally, we implemented the simplification of Section~\ref{sec:main-simplification} for $H\vcentcolon=H_0\cap\SAut(F_4)$, which decreased the number of $H'$-orbits to consider and further increased computation speed by a factor of about $20$ for the problem instances that were most relevant to us.
An annotated but ad-hoc version of the computer code we used to efficiently create the simplified SDP, solve it and verify the solution can be found in the author's GitHub repository.\footnote{\url{https://github.com/MartinNitsche/AutF4-Property-T}}

\medskip

With a well-chosen supporting set $T$, we found a sufficient solution to the simplified SDP for $\SAut(F_4)$, that is, we found -- up to a small numerical error -- a decomposition
\[\Delta^2-\varepsilon\Delta\equiv\sum_{i=1}^n \omega_i^*\omega_i\mod\Span\{\gamma-\gamma h\mid\gamma\in\Gamma,h\in H\},\qquad\text{with }\varepsilon>0.
\]
To conclude that $\SAut(F_4)$ has property~(T) it is now enough to estimate the error, similarly to \cite{NT}, and appeal to Lemma~\ref{lem:main-simplification}. Alternatively, one can use Remark~\ref{rem:primal-simplification} to construct from the solution a decomposition $\Delta^2-\varepsilon'\Delta=\sum_{i=1}^n \hat{\omega}_i^*\hat{\omega}_i$.

In fact, with knowledge of a good supporting set $T$, we solved the corresponding unsimplified SDP, and it turned out that the solution to this much larger SDP works as a witness for property~(T). This solution can be found at \cite{zenodo}. It can be verified in the usual way. To bound the numerical error, it suffices to use the estimate from \cite{NT}. We have attached to this article a short SAGE-script that verifies the solution with interval arithmetic.

\begin{theorem}\label{thm:autF4}
$\Aut(F_4)$, the automorphism group of the free group over four generators, satisfies property~(T).
\end{theorem}

Our implementation of the simplified property~(T) SDP is fast enough that it is possible to experiment with it in a similar way as with $\SL(n,\Z)$. Empirical observations obtained from this might guide the way to a property~(T) proof for $\Aut(F_n)$ that does not rely on the computer.

\begin{question}
In analog to Section~\ref{sec:harpers-operator}, we ask for ``small'' subgroups $\Gamma'<\SAut(F_n)$ such that the SDP property~(T) proof for $\SAut(F_n)$ can be carried out using a support in $\Gamma'$ together with the symmetry of $\SAut(F_n)$.
\end{question}

\section*{Acknowledgments}

This research was supported by ERC Consolidator Grant No.\ 681207 (``Groups, Dynamics, and Approximation'') at TU Dresden and by DFG Grant 281869850 (RTG 2229, ``Asymptotic Invariants and Limits of Groups and Spaces'') at KIT.

The author thanks Andreas Thom for the suggestion of phrasing the dual SDP in the context of Shalom's theorem. This greatly helped to improve the presentation of this article.

\end{document}